\documentclass
{amsart}

\usepackage{eepic}
\newcommand{\dashlinestretch}{10}

\date{February 27, 2000}
\usepackage{a4wide}
\usepackage{amsthm}
\usepackage{amssymb}
\usepackage{mathrsfs}
\swapnumbers
\newcommand{\<}{\langle}
\renewcommand{\>}{\rangle}
\newcommand{\oo}{\text{\it ortho}}

\makeatletter
\newcommand{\KNUTHcases}[1]
{\left \{\,\vcenter {\normalbaselines \m@th
\ialign {$##\hfil $&\quad ##\hfil \crcr #1\crcr }}\right .
}  
\makeatother

\def\newthm#1 {\newtheorem{#1}[thm]{#1}}

\def\itm#1 {\item[{(#1)}]}

\theoremstyle{definition}
\newthm Notation
\newthm {Abuse of Notation}

\newthm {Definition and Fact}
\newthm Definition 
\newthm Fact
\newthm {More Notation}
\newthm {Even More Notation}
\newthm Construction 

\theoremstyle{plain}
\newthm Theorem
\newthm Conclusion
\newthm Lemma
\newthm {Main Lemma}
\newthm Corollary
\newthm Claim

\theoremstyle{remark}
\newthm Observation
\newthm Remark
\newthm Question

\newcommand{\pp}{\perp}

\newcommand{\dom}{{\rm dom}}   
\def\itm#1 {\item[({#1})]}
\def\o0{\mbox{$\{0,1\}$}}

\newcommand{\x}{{\tt x}}
\newcommand{\y}{{\tt y}}
\newbox\happybox
\setbox1=\hbox{$\smile$}
\setbox0=\hbox to \wd1{\hfil$\cdot\cdot$\hfil}

\setbox\happybox=\hbox{\raise4pt\rlap{\box0}\lower1pt\box1}

\newcommand{\xx}{\trianglelefteq}
\newcommand{\dd}{{\it dual}}
\newcommand{\logand}{\wedge}
\newcommand{\F}{{\mathscr F}}

\begin{document}
\author{Martin Goldstern}
\thanks{This research was
supported by the Austrian Science Foundation (FWF), grant P13325-MAT}
\address{Technische Universit\"at\\  Wiedner Hauptstra\3e
  8--10/118.2
\\
A-1040 Wien
}
\email{Martin.Goldstern@tuwien.ac.at}
\urladdr{http://info.tuwien.ac.at/goldstern/}

\title{Interpolation in ortholattices}
 \thanks{This paper is available from {\tt www.arXiv.org}, and also from 
my home page}
\begin{abstract}
If $(L,\vee,\wedge,0,1, {^\pp})$ is a complete ortholattice, $f:L^n\to L$
{\em any} partial function, then there is   
a complete
 ortholattice $L^*$ containing $L$ as a subortholattice, and a
ortholattice polynomial $p$ with coefficients in $L^*$ such that 
$p(a_1,\ldots, a_n) = f(a_1,\ldots, a_n)$ for all $a_1,\ldots , a_n \in
L$.

   Iterating this construction long enough
 yields a complete ortholattice  in which every function can be
 interpolated by a polynomial on any  set  of  small
 enough cardinality.    
 
\end{abstract}

\maketitle




\setcounter{section}{-1}
\section{Introduction}

In \cite{mfl} and \cite{l01} we showed the following:  Let $L$ be any
[bounded] lattice, then there is a lattice $\bar L$ extending $L$
[with the same least and greatest element] such that  every {\it
monotone} function from $L^n$ to $L$ is represented by a polynomial
with coefficients in $\bar L$. 

It is clear that as long as we restrict ourselves to lattice polynomials we
can only interpolate monotone functions.   Here we consider the
problem of interpolation on {\it ortholattices}, i.e., bounded
lattices equipped with an ``orthocomplement''. Since the
orthocomplement reverses order, there is no obvious monotonicity
property that all orthopolynomials in an ortholattice will share. 

The main theorem of this paper shows 
that indeed  there are  no restrictions on the
behavior of orthopolyomials; more precisely:  If $L$ is an
ortholattice, then {\it any} function $f:L^n \to
L$ can be represented by a polynomial with coefficients in some
suitable orthoextension $\bar L$. 

By iterating the construction from the theorem we get, for every
cardinal number $\kappa$,   a lattice $\hat
L$ with the property that  every function from $\hat L^n$ to $\hat L$
can be    
interpolated on any set of size $\le \kappa$.   

We also show that we can construct $\hat L$ such that  $\hat L$ will
be complete (as a partial order).   Moreover, assuming
that the original ortholattice $L$ is complete,
 we construct  $\hat L$ such that  $L$ is a ``convex''
sublattice of $\hat L$.

\section{Basic definitions}

\begin{Notation}
 Lattices are denoted by $L$, $L'$, $L_1$, etc.
When we consider several lattices, we use the self-explanatory
notation $ \logand_{L_1}$ or $\logand_1$,  $\le_2$, etc. for the
operations/relation in 
$L_1$, $L_2$, etc.  
  We agree that the symbol $\wedge$ binds more tightly than $\vee$,
i.e.,  $a\vee b \logand c = a \vee ( b \logand c)$. 

  
For any lattice $L$ we let $L^\dd$ be the dual lattice (with the same
underlying set):   $x \le^{L^\dd} y $ iff $x \ge^L y$.

An ortholattice is a bounded lattice $(L,\vee,\wedge,0,1) $ with an
additional unary  operation $x \mapsto x^\pp$ which satisfies
$x \le y \Rightarrow x^\pp \ge y^\pp$, $(x^\pp)^\pp = x$, 
 $x \vee x^\pp = 1$, $x \wedge x^\pp=0$ for all $x,y$
 (and hence also the de Morgan
laws $(x \vee y)^\pp = x^\pp\wedge y^\pp$, etc.).

\end{Notation}
\begin{Abuse of Notation}
If $(L,\vee,\wedge, 0,1)$ is a bounded lattice, $A \subseteq L$, then 
we say that $A$ is ``convex'' in $L$ iff
\begin{quote}
  whenever $a,a'\in A$, $x\in L$, $0<a\le x\le a'<1$, then also $x\in A$
\end{quote}
i.e., if $A \setminus \{0,1\}$ is convex in $L$ in the usual sense. 

If $L_0 \le L_1$, we say that  $L_0$ is ``downward closed'' in $L_1$
iff:  For all $z\in
      L_0\setminus\{1\}$,  for all $x\in L_1$, if $x\le z $ then $x \in
      L_1$, i.e., if $L_0 \setminus \{1\}$ is downward closed in $L_1$
in the traditional sense. 
\end{Abuse of Notation}

\begin{Definition} Let $L_0$, $L_1$ be bounded lattices. 
\begin{enumerate}
\item $L_0 \le L_1$ means that $L_0$ is a \o0-sublattice of $L_1$
(i.e., $L_0  $ is a sublattice of $L_1 $ with the same maximal and
minimal element)
\item $L_0 \xx L_1 $ means that
 \begin{itemize}
	\itm a $L_0 \le L_1 $
\itm b For every $x\in L_1$ the set $\{z\in L_0: z \le_{L_1} x\}$ has greatest
      element (the ``projection'' of $x$ to $L_0$, written
 $\pi^{L_1}_{L_0}(x)$ or $\pi^1_0(x)$) 
\itm c  $L_0$ is downward closed in $L_1$
 \end{itemize}
\item $L_0 \xx^\dd L_1$ is the dual notion, i.e.,  $L_0^\dd \xx L_1^\dd$.
\end{enumerate}
\end{Definition}

\section{Basic facts}

In this section we collect a few easy facts for later reference.  
We also quote a theorem on interpolation of monotone functions in lattices. 

\begin{Fact}\label{Fjoin}
  Assume $L_0 \xx L_1$, $L_0 \xx^\dd L_2$, $L_1\cap L_2 =
L_0 $.     Let $L= L_1 \cup L_2 $, and let $\le_L  $ be the transitive
closure of $(\le_1) \cup (\le_2) $.  Then
\begin{itemize}
\item $x \le_L y $ iff: $x \le_1 y  $ or $x \le_2 y $ or there exists a
$z\in L_0 $ with $x\le_2 z \le_1 y $. 
\item $(L, \le_L) $ is a lattice with 
$$ 
x \wedge_L y \ = \ 
\KNUTHcases{ x \wedge_1 y              & if $x,y\in L_1 $\cr
              x \wedge_2 y              & if $x,y\in L_2 $\cr
             \pi^1_0(x) \wedge_2 y    & if $x\in L_1 $, $y\in L_2 $\cr
}$$ and similarly for $\vee_L $. 
\end{itemize}
\end{Fact}
\begin{Fact}  If $L_1 \xx L_2 \xx L_3 $, then $L_1 \xx L_3 $.
\end{Fact}

\begin{Fact}\label{equal}
  $L_1$ and $L_2$ are complete lattices,  $L_1 \le L_2$, 
and $A \subseteq L_1$, then $\sup_{L_1} A \ge \sup_{L_2} A$.  

However, if $L_1$ is convex in $L_2$ and $\sup_{L_1} A < 1$,  then 
 $\sup_{L_1} A  =  \sup_{L_2} A$.  
\end{Fact}

\begin{Fact}\label{Flimit}
 Let $(I,\le) $ be a linearly ordered set, and assume that
$(L_i:i\in I) $ is a family of complete lattices such that: 
\begin{itemize}
\item for all  $i<j$, then  $L_i$ is a  $\{0,1\}$-sublattice 
of $L_j$ 
\item for all $i<j$:   $L_i$ is convex in $L_j$. 
\end{itemize}
Then $L:= \bigcup_i L_i  $ is a complete lattice, and $L_i$ is a
convex $\{0,1\}$-sublattice of $   L  $  
for all $  i \in I $. 

\end{Fact}

\begin{proof}

  It is clear that $L $ is a lattice and that 
 $L_i $ is convex in $L$. 

We now check that $L$ is complete.  Let $A \subseteq L$.   We will
show that $\sup_L A$ exists. Wlog we may assume that
$A\not=\emptyset$, $1\notin A$, $0 \notin A$.   Let $a_0\in A\setminus
\{0\}$. 

We may also assume that $1$ is not the least upper bound of $A$, so
let $c<1$ be some upper bound. 

Fix $i_0 $ such that   $a_0, c\in L_{i_0}$. 

We will write $\sup_i$ for the supremum operation in $L_i$.  


Let $i_0 \le i \le j$. Then 
 $ \sup_{i} (A\cap L_i) \le c < 1$, so by fact~\ref{equal} we have 
 $ \sup_{i} (A\cap L_i)
=  \sup_{j} (A\cap L_i)$. 
Hence the sequence $(b_i: i\ge i_1)$, defined by 
$$ b_i:= \sup_i  (A\cap L_i)$$
is weakly increasing. 

For all $i\ge i_1$ we have 
$$ 0 < a_0  \le b_i\le c < 1, \qquad   \qquad  
	a_0,c\in L_{i_0}$$
so since $L_{i_0}$ is convex we get: $\forall i\ge i_0: b_i\in
L_{i_0}$.    
%

Let $a:= \sup_{i_0}(b_i: i\ge i_1$.  Clearly, $a$ is the least upper bound
for $A$.

(Let $a'$ be any upper bound for $A$, say $a'\in L_i$, $i\ge
i_0$, then $a'\ge b_j$ for all $j\ge i$.) 

\end{proof}

\begin{Fact} Let $L = L_1 \cup L_2 $, where 
$L_1 \le  L$, $L_2 \le  L$.   If  $L_1 $ and $L_2 $ are both complete
lattices, then also $L $ is a complete
lattice. 
\end{Fact}

\begin{Theorem}	\label{Told}
 Let $L $ be a complete lattice, $f:L \to L $ a partial
monotone function.

Then there is a complete lattice $\bar L $ and a polynomial $p({\tt
x})\in \bar L[{\tt x}] $ with $p(a)=f(a) $ for all $a\in \dom(f) $. 
Moreover, we can choose $\bar L $ such that 
\begin{itemize}
\item
 $L \xx \bar L $. 
\item $\bar L$ is complete. 
\end{itemize}

\end{Theorem}

\begin{proof} See \cite{l01} and \cite{mfl}.
  (The ``moreover'' part is not stated
there, but following the proof it is easy to see that the lattice
$\bar L $ constructed in \cite{l01} satisfies 
$L \xx \bar L $ and will be complete.)
\end{proof}

\section{From lattices to ortholattices}

We describe a construction that allows us to extend a lattice $L$ to an
ortholattice, preserving the ortholattice structure of a given sublattice 
of $L$.

\begin{Construction}
Let $({L_0},\vee,\wedge,0,1,{^\pp})$ be an ortholattice, $(L_1,
\vee,\wedge,0,1) $ a bounded lattice with ${L_0} \xx L_1$. 

We define two  partial orders
 $L_2 = {\it dual}(L_1, L_0)\supseteq L_0$
 and $L = \oo(L_1, L_0) = L_1 \cup L_2$ as follows: 

Pick a set $L_2$ and a map ${\iota}$ satisfying the following: 
\begin{itemize}
\item ${\iota}:L_1 \to L_2$ is a bijection.
\item $L_1 \cap L_2 = L_0$.
\item ${\iota}(z)=z^{\pp_0} $ for all $z\in L_0$. 
\end{itemize}
(I.e., $L_2\setminus L_0$ is just a disjoint copy of $L_1\setminus
L_0$. 
We make $L_2$ into a lattice by  requiring ${\iota}(x) \le_2 {\iota}(y) $ iff $y
      \le_1 x$, so that  ${\iota}:L_1
      \to L_2$ is a dual isomorphism. 
\end{Construction}

\begin{figure}

\setlength{\unitlength}{0.00083333in}
\begingroup\makeatletter\ifx\SetFigFont\undefined
\def\x#1#2#3#4#5#6#7\relax{\def\x{#1#2#3#4#5#6}}%
\expandafter\x\fmtname xxxxxx\relax \def\y{splain}%
\ifx\x\y   
\gdef\SetFigFont#1#2#3{%
  \ifnum #1<17\tiny\else \ifnum #1<20\small\else
  \ifnum #1<24\normalsize\else \ifnum #1<29\large\else
  \ifnum #1<34\Large\else \ifnum #1<41\LARGE\else
     \huge\fi\fi\fi\fi\fi\fi
  \csname #3\endcsname}%
\else
\gdef\SetFigFont#1#2#3{\begingroup
  \count@#1\relax \ifnum 25<\count@\count@25\fi
  \def\x{\endgroup\@setsize\SetFigFont{#2pt}}%
  \expandafter\x
    \csname \romannumeral\the\count@ pt\expandafter\endcsname
    \csname @\romannumeral\the\count@ pt\endcsname
  \csname #3\endcsname}%
\fi
\fi\endgroup
{\renewcommand{\dashlinestretch}{30}
\begin{picture}(5084,4341)(0,-10)
\put(1936,2069){\ellipse{1800}{4124}}
\path(1425,3756)(150,3756)(1125,2931)
\path(2455,381)(3750,381)(2755,1206)
\put(3300,831){\makebox(0,0)[lb]{\smash{{{\SetFigFont{12}{14.4}{rm}$L_2
= dual(L_1,L_0)$}}}}}
\put(3300,2500){\makebox(0,0)[lb]{\smash{{{\SetFigFont{12}{14.4}{rm}$L
= ortho(L_1,L_0)=L_1\cup L_2$}}}}} 
\put(2100,4131){\makebox(0,0)[lb]{\smash{{{\SetFigFont{12}{14.4}{rm}$
L_0=L_1\cap L_2
$}}}}}
\put(0,3306){\makebox(0,0)[lb]{\smash{{{\SetFigFont{12}{14.4}{rm}$L_1$}}}}}
\end{picture}
}
\end{figure}

      Note that $\le_{0}$ coincides with the
      restriction of $\le_2$ to $L_0$, so $L_0 \le L_2$. 
\relax From  $L_0 \xx L_1$ we conclude $L_0 \xx^\dd L_2$. 

We let 
$L=\oo(L_1, L_0)$ be the (set-theoretic) union of $L_1\cup
L_2$. By fact~\ref{Fjoin} we see that we can make  $L$ into a lattice
 containing $L_1$
and $L_2$ as sublattices.

\begin{Lemma} If $L_0 \xx L_1$, $(L_0, \vee, \wedge, 0,1,{^\pp})$ is
an ortholattice, then there is an ortholattice $L$  such that  
$L_0 \le L_1 \le L$, and $L_0$ is a subortholattice of $L$. 

Moreover, if $L_1$ is complete, then also $L$ will be complete. 
\end{Lemma}

\begin{proof}   We let $L_2 $, ${\iota}$,  and 
$L= \oo(L_1, L_0)$ as above.  

Let $x^{\pp_L} = {\iota}(x) $ for $x\in L_1$ and 
 $ = {\iota}^{-1}(x)$ for $x\in L_2$.

It remains to show that the map
$x\mapsto x^{\pp_L}$ is an orthocomplement for $L$.   Clearly this map
is well defined and an involution, and it agrees with the 
map $x \mapsto x^{\pp_0}$ on $L_0$.  Also, we have $z \vee_L z^{\pp_L} =
z\vee_0 z^{\pp_0} = 1$ for all $z\in L_0$. 

Now let $x\notin L_0$, wlog $x\in L_1$.  We will only check $x \vee_L
x^{\pp_L} = 1$, leaving the dual to the reader. So let $y \ge_L x$, $y
\ge_L  x^{\pp_L}  = {\iota}(x)$. Then $y$ must be in $L_1$, and there is a $z\in
L_0$ such that  ${\iota}(x) \le_2 z \le_1 y$. Now ${\iota}(x) \le_2 z = {\iota}(z^{\pp_0})$ 
implies
$z^{\pp_0} \le_1 x$, hence $ y \ge_1  z \vee  z^{\pp_0} = 1$. 

\end{proof}

\begin{Fact} The operation $\oo({\cdot}, L_0)$ commutes with direct
limits.  In particular, if $(I,{\le})$ is a linear order, $(L_i:i\in
I)$ an increasing family of lattices, $L\xx L_i$ for all $i$, then 
$$ \oo(\bigcup_{i\in I} L_i, L_0) = 
 \bigcup_{i\in I} \oo(L_i, L_0) 
$$
\end{Fact}

\section{Theorems}

We prove the two main theorems mentioned in the introduction. 
We conclude with an open question concerning the difference of unary
and $n$-ary functions. 

\begin{Theorem}\label{extend}
If $L_0$ is an ortholattice, $f:L_0 \to L_0$, then there is an ortholattice
$L^*$ extending $L_0$ such that  $f$ is the restriction of a polynomial
function over $L^*$.

Moreover, if $L_0$ is complete then we can have $L^* = \oo(L_1, L_0)$ for
some complete $L_1$, $L_0 \xx L_1$. 

\end{Theorem}

\begin{proof}
Since every ortholattice can be embedded into a complete ortholattice 
(the MacNeill completion; see, e.g., \cite[4.1]{Bruns:1976},
\cite{Kalmbach:1983})
we may assume that $L_0$ is complete. 
Let $f:L_0 \to L_0$. 

Let $L_0'$ be horizontal sum of $L_0 $ and $L_0 \times L_0$, i.e.,
assume that $L_0$ and $L_0 \times L_0$ have the same least and greatest
elements (but are otherwise disjoint), and make $L_0' = L_0 \cup (L_0
\times L_0)$ into a lattice by taking 
$$\le_{L_0'} = (\le_{L_0}) \cup (\le_{L_0\times L_0}) .$$
Note that $L_0'$ is a complete \o0-lattice and $L_0 \xx L_0'$. 

Now consider the partial functions $\bar f$, $g_1$ and $g_2$, defined by 
\begin{itemize}
\item $\bar f( \<x,x^\pp\>) = f(x)$  for all $x\in L_0$. 
\item $g_1(x) =  \<x,0\>$, $g_2(x) =  \<0,x\>$, 
 for all $x\in L_0$. 
\end{itemize}

Notice the elements of the set $\{ \<x,x'\> : x \in L_0\}$ are
pairwise incomparable, so the function $\bar f$ is trivially monotone.

By theorem~\ref{Told} we can find a lattice $L_1$, $L_0' \xx L_1$ in
which the functions $\bar f$, $g_1$  and $g_2$ are restrictions of polynomials
$p$, $q_1$  and $q_2$, respectively.  
Now let $L = \oo(L_1, L_0)$, so  $L$  is
an orthoextension of $L_0$. 

Now $h({\x}) = p(q_1(\x) \vee q_2(\x^\pp)) $ is an orthopolynomial
with coefficients in $L$, and clearly    $h(x) = p ( 
\<x,0\> \vee \<0,x^\pp\>) =
p(\<x,x^\pp\>) =
\bar f(\<x,x^\pp\>) =  f(x)$ for all $x\in L_0$. 

\end{proof}

\begin{Remark}
For every orthopolynomial $p(\x)$ there is a lattice polynomial
$p'(\x,\y)$ such that  (by de Morgan's laws) $p(\x)$ is equivalent to
$p'(\x,\x^\perp) $.
\end{Remark}

\begin{Definition} \label{power}
Let  $\F$ a family of lattices.   We say that 
$\F$ is power closed if: 
\begin{quote}
  For every $S\in \F$ there is some $S'\in \F$ which is isomorphic to
  $S\times S$. 
\end{quote}
We say that $L$ is $\kappa$-power closed if the family of sublattices
of size $\le \kappa $ is power closed. 
\end{Definition}
\begin{Fact}
If $\F$ is power closed, then: 
  For every $S\in \F$ and every $n>1$ there is some $S'\in \F$ which
  is isomorphic to   $S^n$. 
\end{Fact}

\begin{Theorem}
Let $L$ be a complete ortholattice, and let $\kappa$ be any cardinal.
 Then:
\begin{enumerate} 
\item
  There is a complete ortholattice
$\bar L$ extending $L$ such that  
 every function
$f:L \to L$ is represented by an orthopolynomial of $\bar L$. 
\\
Moreover, $\bar L$ can be chosen to be of the form $\oo
(L_1,L)$ with $L \xx L_1$.   In particular, $L$ will be convex in $\bar
L$. 
\\
Moreover, $\bar L$ can be chosen to be $\kappa$-power closed. 
\item There is a complete
ortholattice $\hat L$ extending $L$ such that: 
\begin{quote}
For every natural number $n$,  for every  function 
$f: \hat L^n \to \hat L$  and for every set $A \subseteq L^n$ of
cardinality $\le \kappa$ there is an orthopolynomial $p({\tt x}_1, 
\ldots , {\tt x}_n)$ with coefficients in $\bar L$ that interpolates
$f$ on every point in $A$. 
\end{quote}
\end{enumerate}
\end{Theorem}

\begin{proof}
Choose a cardinal $\lambda$ of cofinality $>\kappa $ such that 
there is a transfinite  enumeration (not necesasrily 1-1)
$\{f_i: 0 \le i < \lambda\}$ of all 
functions $f:L \to L$. 
  Define an increasing
transfinite sequence $(L_i: i \le  \lambda)$ of $\{0,1\}$-lattices 
satisfying 
\begin{enumerate}
\item $L_0 = L$. 
\item If $i< j\le \lambda$,  then $L_i \xx L_j$
\item For every $i< \lambda$ there is a lattice polynomial $p_i$ with
coefficients in $L_{i+1}$ such that  for all $z\in L$: $f_i(z) = p_i(z,
z^{\pp_L})$. 
\item For every $i< \lambda$, $L_{i+1}$ contains an isomorphic copy 
of $L_i^2$. 
\item If $i$ is a limit stage, then $L_i$ is the direct limit of
$(L_j:j<i)$,  (i.e., $L_i = \bigcup _{j<i} L_j$).
\end{enumerate}

Finally, let $\bar L = \oo(L_\lambda, L)$. 
Note that $\bar
L$ will contain an isomorphic copy of $\oo(L_i,L)$ for every $i$.
This finishes the proof
of the first statement.   

To prove the second claim, apply the conclusion from the first claim 
$\kappa^+$ many times to get the conclusion for all unary functions, 
 and then  use fact \ref{n} to take care of all $n$-ary functions. 
\end{proof}

\begin{Fact}\label{n}
Fix a cardinal number $\kappa$.  We call a set $A$ ``small'' if the
cardinality of $A$ is $<\kappa $. 

Assume that $L$ $\kappa$-power closed, and that (a) or (b) holds: 
\begin{itemize}
\itm a $L$ is a lattice, and for every small $A \subseteq L$ and for 
every {\em monotone}  $f:A \to L$ there is  a lattice polynomial
$p\in L[{\tt x}]$ such that  $f(a)=p(a)$ for all $a\in A$
\itm b $L$ is an ortholattice, and for every small $A \subseteq L$ and for 
{\em every}  $f:A \to L$ there is  an orthopolynomial
$p\in L[{\tt x}]$ such that  $f(a)=p(a)$ for all $a\in A$
\end{itemize}
Then also (a') or (b'), respectively, holds: 
\begin{itemize}
\itm a'  
For every $n$, for every small $A
 \subseteq L^n$,
 for every monotone $g:A \to L$
 there 
 is a lattice polynomial
$p\in L[{\tt x_1}, \ldots, {\tt x}_n]$
 such that  $$g(a_1, \ldots a_n)=p(a_1, \ldots, a_n)
\mbox{ for all }(a_1, \ldots , a_n)\in A.$$
\itm b'  Like (a'), but again replace ``every monotone'' by ``every'',
and 
``lattice polynomial'' by ``orthopolynomial''. 
\end{itemize}
\end{Fact}

\begin{proof} Let $g:L^n\to L$, $A \subseteq L$ small.  Let $S$ be the
 sub(ortho)lattice generated by $A$.  Since $L$ is $\kappa$-power closed, we
 can find a sublattice $S' \le L$ and an isomorphism $\iota:S^n\to S'$. 
For $\ell=1,\ldots, n$ let $\iota_\ell:S\to S'$ be defined by 
$$ \iota_\ell(s) = \iota(0,\ldots, 0, s, 0, \ldots, 0)
\qquad \mbox{here, $s$ appears in the $\ell$-th coordinate}$$
Note that
 $$(*) \qquad \qquad 
\iota(s_1,\ldots, s_n) = \iota_1(s_1)\vee \cdots \vee
\iota_n(s_1)$$ 

Define  a partial function $f:L\to L$ as follows: 
$$(**) \qquad f(i(s_1, \ldots, s_n)) = g(s_1, \ldots, s_n)$$
and  $f(t) = $ undefined if $t\notin S'$.   Note that $f$ is monotone
if $g$ is monotone.

Now note that $g$, $\iota_1$, \dots, $\iota_n$ are unary functions
from $L$ to $L$, so they are all represented by polynomials, and
so by $(*)$ and $(**)$ also $g$ is represented by a polynomial. 
\end{proof}

Note that the construction we used in theorem \ref{extend} will
automatically guarantee that the resulting structure will be 
$\kappa$-power closed (see definition~\ref{power}), so by fact \ref{n}
there will be no difference between unary and $n$-ary interpolation. 
This motivates the following questions: 

\begin{Question}
\begin{enumerate}
\item Are there infinite ortholattices $L$ where every function
$f:L\to L  $ can be interpolated by an orthopolynomial 
on every (say) countable set, 
but not every function $f:L^2 \to L$? 
\item Are there infinite lattices $L$ where every monotone function
$f:L\to L  $ can be interpolated by a lattice polynomial 
on every (say) countable set, 
but not every monotone function $f:L^2 \to L$? 
\end{enumerate}
\end{Question}

\nocite{gratzer:1998}


\end{document}